\documentclass[12pt, a4paper]{article}
\usepackage[T1]{fontenc}
\usepackage[french]{babel}
\usepackage{amsmath,amsthm}
\usepackage{amssymb}
\usepackage{amscd}
\usepackage[latin1]{inputenc}
\usepackage{amsfonts}
\usepackage[matrix,arrow,curve]{xy}

\theoremstyle{plain}
\newtheorem{theo}{Th\'eor\`eme}[section]
\newtheorem{lem}[theo]{Lemme}
\newtheorem{prop}[theo]{Proposition}
\newtheorem{cor}[theo]{Corollaire}

\theoremstyle{definition}

\newtheorem{theosd}[theo]{Th\'eor\`eme}
\newtheorem{propsd}[theo]{Proposition}
\newtheorem{rem}[theo]{Remarque}

\newcommand{\Zn}{\ensuremath{\mathbb Z/n\mathbb Z}}

\newcommand{\Zl}{\ensuremath{\mathbb Z/l\mathbb Z}}
\newcommand{\Zln}{\ensuremath{\mathbb Z/l^n\mathbb Z}}
\newcommand{\FF}{\ensuremath{\mathbb F}}

\newcommand{\Hn}{\ensuremath{H^3_{\mathrm{nr}}}}
\newcommand{\nH}{\ensuremath{H_{\mathrm{nr}}}}

\newcommand{\QZl}{\ensuremath{\mathbb Q_l/\mathbb Z_l}}
\newcommand{\QZq}{\ensuremath{\mathbb Q_{l'}/\mathbb Z_{l'}}}

\newcommand{\eqdef}{\ensuremath{\stackrel{\mathrm{d\acute{e}f}}{=}}}

\title{Cohomologie non ramifiée en degré trois d'une variété de Severi-Brauer}
\author{Alena Pirutka}

\begin{document}

\maketitle

\begin{abstract}
 Soit $K$ le corps des fractions d'une surface projective et lisse, géométriquement intègre, définie sur un corps fini $\FF$, $\mathrm{car}.\FF\neq 2$. Soit $C/K$ une conique. Parimala et Suresh \cite{PS}  ont montré que le groupe de cohomologie non ramifiée $\Hn(K(C)/\FF, \QZl(2))$ est nul pour tout $l\neq{\mathrm{car}.\FF}$.  Dans cette note on étend leur résultat aux variétés de Severi-Brauer associées à une algèbre centrale simple dont d'indice $l$ est premier et différent de $\mathrm{car}.\FF$.
\end{abstract}

\section{\large{Introduction}}

Soit $K$ un corps.  Pour $n$ un entier  inversible sur $K$,  on note $\mu_{n}$ le $K$-sch\'ema en groupes (\'etale) des racines $n$-i\`emes de l'unit\'e. Pour $j$ un entier positif on note $\mu_{n}^{\otimes j}=\mu_{n}\otimes\ldots\otimes\mu_{n}$ ($j$ fois). On pose $\mu_{n}^{\otimes j}=Hom_{k-gr}(\mu_{n}^{\otimes (-j)}, \mathbb Z/n)$ si $j$ est n\'egatif et $\mu_{n}^{\otimes 0}=\mathbb Z/n$.  Ces $k$-sch\'emas en groupes donnent des faisceaux \'etales, not\'es encore $\mu_{n}^{\otimes j}$, sur toute $K$-vari\'et\'e $X$. On note $H^i(X,\mu_n^{\otimes j})$ les groupes de cohomologie \'etale de $X$ \`a valeurs dans $\mu_n^{\otimes j}$. Lorsque $K$ contient une racine primitive $n$-ième de l'unité, on a un isomorphisme   $\mu_{n}^{\otimes j}\stackrel{\sim}{\to}\Zn$ pour tout $j$.

Soit $K$ un corps. Soit $F$ un corps de fonctions sur $K$.  Soient $j\geq 1$ un entier naturel  et $i\in \mathbb Z$ un entier relatif. Dans la suite on va utiliser les notions suivantes de  cohomologie \textit{non ramifi\'ee} :
\begin{itemize}
 \item[(i)] $\nH^j(F/K, \mu_n^{\otimes i})=\bigcap\limits_A\mathrm{Ker}[H^j(F, \mu_n^{\otimes i})\stackrel{\partial_{j,A}}{\to}H^{j-1}(k_A, \mu_n^{\otimes i-1})].$
Dans cette formule, $A$ parcourt les anneaux de valuation discr\`ete de rang un, de  corps des fractions $F$, contenant le corps $K$.  Le corps r\'esiduel d'un tel anneau $A$ est not\'e $k_A$ et l'application $\partial_{j,A}$ est l'application r\'esidu.
\item [(ii)] Pour $X$ une $K$-vari\'et\'e int\`egre, on note $\nH^j(X, \mu_n^{\otimes i})\eqdef\nH^j(K(X)/K, \mu_n^{\otimes i})$.

Pour $X$  propre et lisse, les r\'esultats de Bloch et Ogus permettent d'identifier  $\nH^j(X, \mu_n^{\otimes i})$ au groupe de cohomologie de Zariski $H^0(X, \mathcal H^j(\mu_n^{\otimes i}))$, où $\mathcal H^j(\mu_n^{\otimes i})$ d\'esigne le faisceau de Zariski sur $X$ associ\'e au pr\'efaisceau $U\mapsto H^j(U,\mu_n^{\otimes i})$ (cf. \cite{CT}).
\item[(iii)] Si $X$ est régulier en codimension $1$, on pose $$\nH^j(K(X)/X, \mu_n^{\otimes i})=\bigcap\limits_{x\in X^{(1)}}\mathrm{Ker}\, \partial_{j,\mathcal{O}_{X,x}}$$    la cohomologie non ramifiée de $K(X)$ par rapport à $X$. \\
\end{itemize}


Soit maintenant $K$ le corps des fractions d'une surface lisse géométriquement intègre $S$, définie sur un corps fini $\FF$. Dans \cite{PS}, Parimala et Suresh montrèrent que $\Hn(K(X)/\FF, \QZl(2))=0$, $l\neq{\mathrm{car}.\FF}$,  pour $X$ une conique sur $K$.  Le but de cette note est d'étendre leurs arguments au cas des variétés de Severi-Brauer d'indice premier :

\theosd\label{main}{\textit{Soit $K$ le corps des fractions d'une surface projective et lisse, géométriquement intègre $S$, définie sur un corps fini $\FF$. Soit $X$ la variété de Severi-Brauer associée à un corps gauche de centre $K$ et d'indice premier  $l\neq{\mathrm{car}.\FF}$. On a alors  $$\Hn(K(X)/\FF, \QZq(2))=0$$ pour tout premier $l'$,  $l'\neq{\mathrm{car}.\FF}$.\\}}

\rem{Ce résultat est aussi vrai pour une variété de Severi-Brauer associée à une $K$-algèbre centrale simple d'indice premier  $l\neq{\mathrm{car}.\FF}$. En effet,  $\Hn(X, \QZq(2))=\Hn(X\times \mathbb P^n_K, \QZq(2))$ d'après \cite{CTOj} 1.2.\\ }

Pour montrer ce théorème, on doit essentiellement établir que  $\Hn(K(X)/\FF, \mu_l^{\otimes 2})=0$ (cf. section $2$). Pour ce faire, on suit les m\^emes étapes que dans \cite{PS}. On montre d'abord que tout élément  $\beta\in\Hn(K(X)/\FF, \mu_l^{\otimes 2})$ provient d'un élément    $\xi\in H^3(K, \mu_l^{\otimes 2})$ (cf. section $2$). Dans la section $3$, en utilisant le principe de type local-global de \cite{PS}, on montre que l'on peut en effet supposer que $\xi\in H^3_{nr}(K/S, \mu_l^{\otimes 2})$. D'après le théorème de Colliot-Thélène, Sansuc et Soulé \cite{CTSS}, le groupe $H^3_{nr}(K/S, \mu_l^{\otimes 2})$ est nul, ce qui nous permet de conclure (cf. section $3$).  Dans la section $4$, on donne aussi quelques détails de la preuve du principe local-global de \cite{PS}. \\

\section{\large{Comparaison entre cohomologie non-ramifiée en degré trois d'une variété de Severi-Brauer et cohomologie du corps de base}}



\prop\label{comp1}{Soit $K$ un corps. Soit $X$ la variété de Severi-Brauer associée à un corps gauche $D$ de centre $K$ et d'indice premier  $l\neq{\mathrm{car}.K}$.
L'application naturelle $$H^3(K, \QZl(2))\stackrel{\phi}{\to} \Hn(X, \QZl(2))$$ est surjective.  }
\rem{D'après un résultat de Peyre \cite{Pe}, on comprend explicitement le noyau de $\phi$ : un élément de $\mathrm{ker} \phi$ s'écrit comme un cup-produit $\alpha\cup f$ où $\alpha\in Br\,K$ est la classe de $D$ et $f\in K^*$. Cela ne nous est pas nécessaire pour la suite.}
\proof{Cette proposition est une conséquence des résultats de B. Kahn \cite{K1}, \cite{K2}. Supposons d'abord que $K$ est parfait. Soient $K^s$ une clôture séparable de $K$, $\bar X=X\times_K K^s$ et $G=Gal(K^s/K)$. D'après \cite{K1} 5.3(8) et \cite{K2} 2.5, pour $X$ une $K$-variété de Severi-Brauer, on a alors un complexe :
$$0\to \mathrm{Coker} \phi\,\{l\}\to \mathrm{Coker}[CH^2(X)\to CH^2(\bar X)^G]\{l\}\stackrel{\delta}{\to} Br\,K\{l\}, $$
qui est exact sauf peut-être en terme au milieu. Si $X$ est une conique, on a immédiatement $\mathrm{Coker} \phi\,\{l\}=0$ car  $CH^2(\bar X)=0$ (ce résultat est du à Suslin, \cite{Su}). 

Supposons que $l\geq 3$. Puisque $D$ est d'indice premier $l$, on a que $\mathrm{Coker}[CH^2(X)\to CH^2(\bar X)^G]\simeq \Zl$ (\cite{MS}, 8.7.2). Dans ce cas, d'après \cite{K1} 7.1 et \cite{K2} 2.5, on a explicitement $\delta(1)=2[D]\neq 0$, où $[D]$ désigne la classe de $D$ dans $Br K$ et  $1\in \Zl$ est identifié au générateur de $\mathrm{Coker}[CH^2(X)\to CH^2(\bar X)^G]$. Ainsi, sur un corps $K$ parfait, l'application $H^3(K, \QZl(2))\stackrel{\phi}{\to} \Hn(X, \QZl(2))$ est surjective.

Dans le cas général, soit $K'$  une clôture parfaite de $K$. On a le diagramme commutatif

\begin{center}
\xymatrix{
 H^3(K, \QZl(2))\ar[d]\ar[r]^{\phi} &\Hn(X, \QZl(2))\ar[d]&\\
 H^3(K', \QZl(2))\ar[r]^{\phi'} &\Hn(X_{K'}, \QZl(2)).&
}
\end{center}
Soit $\beta\in \Hn(X, \QZl(2))$. D'après ce qui prècede, l'application $\phi'$ est surjective. Ainsi, il existe une extension finie $K''/K$ dont le degré est une puissance de ${\mathrm{car}.K}$, telle que $Res_{K''/K}(\beta)$ provient d'un élément de $ H^3(K'', \QZl(2))$. Puisque $l\neq{\mathrm{car}.K}$, on obtient le résultat par un argument de  corestriction.

 \qed\\}


\cor\label{comp2}{Soit $K$ le corps des fractions d'une surface projective et lisse, géométriquement intègre $S$, définie sur un corps fini $\FF$. Soit $X$ la $K$-variété de Severi-Brauer aassociée à un corps gauche de centre $K$ et d'indice premier  $l\neq{\mathrm{car}.\FF}$. Alors
\begin{itemize}
\item[(i)] pour tout $l'\neq l,\mathrm{car}.K $, on a  $\Hn(K(X)/\FF, \QZq(2))=0$;
\item[(ii)] $\Hn(K(X)/\FF, \QZl(2))=\Hn(K(X)/\FF,\mu_l^{\otimes 2})$ et l'image de l'application naturelle $$H^3(K, \mu_l^{\otimes 2})\to H^3(K(X), \mu_l^{\otimes 2})$$ contient le groupe  $ \Hn(K(X)/\FF, \mu_l^{\otimes 2})$.
\end{itemize} }
\proof{ Soit $K'$ une extension de $K$ d'indice $l$, telle que $X_{K'}$ est isomorphe à un espace projectif. Par le même argument que dans \cite{CT} 2.1.10, on a une application entre les cohomologies non ramifiées 
\vspace{-0.3cm}
$$\Hn(K(X)/\FF, \QZq(2))\to H^3_{nr}(K'(X_{K'})/\FF,\QZq(2)).$$
Puisque $K'(X_{K'})$ est une extension transcendante pure de $K'$, on a un isomorphisme (\cite{CTOj} 1.2) :
\vspace{-0.4cm}
 $$H^3_{nr}(K'/\FF, \QZq(2))\stackrel{\sim}{\to}\Hn(K'(X_{K'})/\FF, \QZq(2)).$$
 Le groupe $H^3_{nr}(K'/\FF, \QZq(2))$ est nul d'après \cite{CTSS} p.790.  Un argument de corestriction montre alors que pour $l'\neq l$, on a $\Hn(K(X)/\FF, \QZq(2))=0$, et que tout élément de $\Hn(K(X)/\FF, \QZl(2))$ est annulé par $l$. On a alors que tout élément de $\Hn(K(X)/\FF, \QZl(2))$ est annulé par $l$ et il vient donc de $\Hn(K(X)/\FF, \mu_l^{\otimes 2})$ par \cite{MS} (p. 339).

Soit $\xi\in\Hn(K(X)/\FF, \mu_l^{\otimes 2}).$ On voit  $\xi$ aussi comme un élément de $\Hn(X,\mu_l^{\otimes 2})$ et on déduit de la proposition précédente que $\xi$ provient d'un élément $\beta\in H^3(K, \QZl(2))$.  Montrons que $\beta$ est annulé par $l$.  


Par le même raisonnement que précédemment, l'image $\xi'$ de $\xi$ dans $H^3(K'(X_{K'}), \mu_l^{\otimes 2})$ est nulle. En effet, 
$\xi'\in \Hn(K'(X_{K'})/\FF, \mu_l^{\otimes 2})\stackrel{\sim}{\to} \Hn(K'/\FF, \mu_l^{\otimes 2})=0$.  On en déduit que l'image de $\beta$ dans $H^3(K', \QZl(2))$ est nulle. On a alors : $l\beta=Cor_{K'/K}\circ Res_{K'/K}(\beta)=0$. D'après \cite{MS}, on a alors que $\beta$ vient d'un élément de  $H^3(K, \mu_l^{\otimes 2})$.

\qed\\}

\section{\large{Preuve du théorème \ref{main}}}
Montrons d'abord le lemme suivant :

\lem\label{degn}{Soit $R$ un anneau de valuation discrète dont on note $K$ le corps des fractions et $k$ le corps résiduel.  Soit $X$ la variété de Severi-Brauer associée à un corps gauche $D$ de centre $K$ et d'indice premier $l$, $(l,\mathrm{car}.\, k)=1$. Soit $\alpha$ la classe de $D$ dans $Br\,K$.  Soit $\xi\in H^3(K,\mu_l^{\otimes 2})$. Supposons que pour toute valuation $v$ sur $K(X)$ induisant sur $K$ soit la valuation triviale, soit la valuation associée à $R$, le résidu $\partial_v(\xi_{K(X)})$ est nul.   On a alors :
\begin{itemize}
 \item [(i)] si $\alpha$ est non ramifiée en $R$, alors $\partial_R(\xi)$ est un multiple de la spécialisation  $\bar{\alpha}=\partial_{R}(\alpha\cup \pi)$,  où l'on note $\pi$ une uniformisante de $R$;
\item[(ii)] si $\alpha$ est ramifiée en $R$, alors ou bien $\xi$ est non ramifiée en $R$, ou bien  $\partial_R(\xi)$ est isomorphe à une algèbre cyclique $(\partial_R(\alpha), c)$ pour $c\in k^*$.
\end{itemize}
}
\proof{  Pour prolonger la valuation sur $K$ en une valuation sur $K(X)$ on s'intéresse à la structure de la fibre spéciale d'un modèle  de $X$ au-dessus de $R$.

Quitte à changer $K$ par son complété, on peut supposer que $K$ est complet.  Notons d'abord que $\alpha$ est trivialisée par une extension non ramifiée de $K$. En effet, soit  $K_{nr}$  l'extension maximale non ramifiée de $K$, soit $R_{nr}$ l'anneau des entiers de $K_{nr}$ et soit $k^s$ une clôture séparable de $k$. On a une suite exacte $$0\to H^2(R_{nr},\mu_{l^n})\to H^2(K_{nr},\mu_{l^n})\stackrel{\partial_R}{\to} H^1(k^s,\Zln).$$ Puisque $k^s$ est séparablement clos, on a $H^1(k^s,\Zln)=0$ et
$ H^2(R_{nr},\mu_{l^n})=H^2(k^s,\mu_{l^n})=0$, d'où $Br\,K_{nr}\{l\}=0$  et donc $\alpha$ est trivialisée par une extension non ramifiée de $K$.

Supposons que $\alpha$ est non ramifiée. Dans ce cas, puisque $\alpha$ est trivialisée par une extension non ramifiée de $K$, on a que  $D$ se prolonge en une algèbre d'Azumaya $\Lambda$ sur $R$, qui a pour la classe $\alpha$ (cf. \cite{Fr}, 1.1 et 1.2). Ainsi $X$ se prolonge en un schéma de Severi-Brauer au-dessus de $R$, associé à $\Lambda$, dont la fibre spéciale $\bar X$ est donnée par l'image de $\alpha$ dans $H^2(k,\mu_{l})$, qui est précisement $\bar{\alpha}$. 

Supposons que $\alpha$ est ramifiée.    Sous l' hypothèse que $\alpha$ est trivialisée par une extension non ramifiée de $K$,  d'après Artin \cite{A} 1.4, il existe un modèle de $X$ au-dessus de $R$ dont la fibre spéciale géométrique est réduite et contient $l$ composantes conjuguées sur $k$, qui sont des variétés rationnelles. 

On voit ainsi qu'il existe une valuation $v$  sur $K(X)$ qui prolonge la valuation sur $K$ dont le corps résiduel $\kappa(v)$ est
\begin{itemize}
 \item [(i)]$\kappa(v)=k(\bar X)$, ou  $\bar X$ est une $k$-variété de Severi-Brauer de classe $\bar{\alpha}$, si $\alpha$ est non ramifiée;
\item [(ii)]$\kappa(v)$ est une extension transcendante pure d'une extension $k'$ de $k$ de degré $l$, sinon.
\end{itemize}
On a alors le diagramme commutatif suivant :
\begin{center}
\xymatrix{
H^3( K, \mu_l^{\otimes 2})\ar[r]^{Res}\ar[d]^{\partial_R} &H^3(K(X),\mu_l^{\otimes 2})\ar[d]^{\partial_v}&\\
H^{2}(k, \mu_l)\ar[r]^{Res} & H^{2}(\kappa(v), \mu_l). &
}
\end{center}
Puisque $\xi$ devient non ramifiée sur $K(X)$, on a que $\partial_R(\xi)$ est dans le noyau de l'application $H^{2}(k, \mu_l)\to H^{2}(\kappa(v), \mu_l)$. Si $\alpha$ est non ramifiée,  $\partial_R(\xi)$ est alors  un multiple de $\bar {\alpha}$ d'après le théorème d'Amitsur, ce qui établit $(i)$.

Supposons maintenant que $\alpha$ et $\xi$ sont ramifiées en $R$.  Puisque  $\kappa(v)$ est une extension transcendante pure de $k'$, $H^{2}(k', \mu_l)$ s'injecte dans $H^{2}(\kappa(v), \mu_l)$. Ainsi $\partial_R(\xi)$ est dans le noyau de l'application $H^{2}(k, \mu_l)\to  H^{2}(k', \mu_l)$.

D'autre part, puisque $\alpha$ devient triviale sur $K(X)$,  on voit que $\partial_R(\alpha)$ est dans le noyau de l'application $H^{1}(k, \Zl)\to  H^{1}(k', \Zl)$ par le même argument. Puisque $\partial_R(\alpha)$ est un élément non nul dans $H^{1}(k, \Zl)$, il correspond à une extension  galoisienne, cyclique, de degré $l$, qui coïncide avec $k'$ car $\partial_R(\alpha)\in \mathrm{Ker}[H^{1}(k, \Zl)\to  H^{1}(k', \Zl)]$ et $[k':k]=l$.  Puisque $\partial_R(\xi)$ est dans le noyau de l'application $H^{2}(k, \mu_l)\to  H^{2}(k', \mu_l)$, cela implique que  $\partial_R(\xi)$ est isomorphe à une algèbre cyclique $(\partial_R(\alpha), c)$ pour $c\in k^*$ (\cite{Se}, p.211), ce qui établit $(ii)$.

  \qed\\}

Passons maintenant à la preuve du théorème \ref{main}. D'après \ref{comp2}, il suffit  d'établir que  $\Hn(K(X)/\FF, \mu_l^{\otimes 2})=0$. Notons qu'on peut supposer que $K$ contient une racine primitive $l$-ième de l'unité. En effet,   le degré $d$ de l'extension $K'$ de $K$, obtenue en ajoutant une racine primitive $l$-ième de l'unité, divise $l-1$. Ainsi $d\mathrm{Id}=Cor_{K'(X)/K(X)}\circ Res_{K'(X)/K(X)}$ est un isomorphisme sur $\Hn(K(X)/\FF, \mu_l^{\otimes 2})$. Il suffit donc d'établir $\Hn(K'(X)/\FF, \mu_l^{\otimes 2})=0$.

Soit $\alpha$  la classe de $D$ dans $Br K$.
Soit $\beta\in \Hn(K(X)/\FF, \Zl)$. D'après \ref{comp2}, $\beta$ provient d'un élément $\xi\in H^3(K, \Zl)$. Montrons qu'il existe $f\in K^*$ tel que $\xi'=\xi-\alpha\cup f\in \Hn(K/S,\Zl)$. Pour ce faire, on utilise le théorème \ref{lg} ci-après (principe local-global de \cite{PS}).
D'après ce théorème, il suffit de trouver, pour tout point $x\in S$ de codimension $1$, un élément  $f_x\in K_x^*$ tel que $\xi-\alpha\cup f_x\in \Hn(K_x,\Zl)$. On a trois cas à considérer :
\begin{enumerate}
\item $\xi$ est non ramifiée en $x$. Dans ce cas,  $f_x=1$ convient.
\item $\xi$ est ramifiée en $x$ et $\alpha$ est non ramifiée en $x$. D'après le lemme \ref{degn}, $\partial_x(\xi)=r\bar \alpha$. Soit $\pi$ une uniformisante de $\mathcal O_{S,x}$. Alors $f_x=\pi^{r}$ convient : $\partial_x(\xi-\alpha\cup \pi^r)=r\bar \alpha-r\bar \alpha=0$.
\item $\xi$ est ramifiée en $x$ et $\alpha$ est ramifiée en $x$. D'après le lemme \ref{degn}, on peut écrire $\partial_x(\xi)=(\partial_x(\alpha), c)$. On relève $c$ en une unité $c'\in K_x$. Puisque  $\partial_x(\alpha\cup c')=(\partial_x(\alpha), c)$,  $f_x=c'$ convient.
\end{enumerate}

Ainsi il existe $f\in K^*$ tel que $\xi'=\xi-\alpha\cup f\in \Hn(K/S,\Zl)$.   On a ainsi que $\beta$ provient de $\xi'\in H^3_{nr}(K/S, \mu_l^{\otimes 2}).$ D'après \cite{CTSS}, p.790, ce groupe est nul (cf. aussi \cite{CT} 2.1.8). Ainsi $\Hn(K(X)/\FF, \mu_l^{\otimes 2})=0$, ce qui termine la preuve du théorème \ref{main}.\qed

\section{\large{Le principe local-global de Parimala et Suresh dans le cas d'indice premier}}

Le théorème  suivant  est démontré dans \cite{PS} seulement dans le cas où $\alpha$ est un symbole. Rappelons ici la preuve pour  nous assurer que l'hypothèse que $\alpha$ est d'indice $l$ suffit.

\theosd\label{lg}{\textit{Soit $K$ le corps des fractions d'une surface projective et lisse, géométriquement intègre $S$, définie sur un corps fini $\FF$. Soit $l$ un entier premier, $(l,\mathrm{car}. K)=1$. Supposons que $K$ contient une racine primitive $l$-ième de l'unité.  Soit $\alpha\in H^2(K,\Zl)$ un élément d'indice $l$ et soit $\xi\in H^3(K, \Zl)$. Les assertions suivantes sont équivalentes :
\begin{itemize}
 \item [(i)] il existe $f\in K^*$ tel que $\xi-\alpha\cup f\in \Hn(K/S,\Zl)$;
\item[(ii)] pour tout point $x\in S$ de codimension $1$, il existe un élément non nul $f_x$ dans le complété $K_x$ de $K$ en $x$ tel que  $\xi-\alpha\cup f_x\in \Hn(K_x,\Zl)$.\\
\end{itemize}}
  }

Soit $K$ le corps des fractions d'une surface régulière propre, géométriquement intègre $S$, définie sur un corps  $k$. Soit $l$ un entier premier, $(l,\mathrm{car}. k)=1$. Supposons que $K$ contient une racine primitive $l$-ième de l'unité.  Soit $\alpha\in H^2(K,\Zl)$ un élément d'indice $l$. Pour démontrer le principe \ref{lg} on va utiliser la théorie de ramification d'algèbres à division, développée par Saltman (\cite{S1}, \cite{S2}).

Soit $ram_S(\alpha)$ le diviseur de ramification de $\alpha$ dans $S$. D'après la résolution des singularités de surfaces, quitte à éclater $S$, on peut supposer que le support de $ram_S(\alpha)$ est l'union des courbes régulières intègres à croisements normaux. On note $C_1,\ldots C_n$ ces courbes. D'après Saltman, après éventuellement quelques éclatements, on peut associer à chaque $C_i$ son coefficient $s_i$. La construction de $s_i$ tient compte de la ramification de $\alpha$ au-dessus des points d'intersections des divers $C_i$ et $C_j$. Pour la suite on n'aura pas besoin de détailler cette construction.

Soient $F_1,\ldots F_r$ des courbes irréductibles régulières et propres dans $S$, telles que $\{C_1,\ldots C_n, F_1,\ldots F_r\}$ soient à croisements normaux. Soient $m_1,\ldots m_r$ des entiers. L'assertion suivante est démontrée dans \cite{PS}, 2.2 (et utilise  \cite{S1}, 4.6 et \cite{S2}, 7.8) :

\propsd\label{div}{\textit{Avec les notations précédentes, il existe $f\in K^*$ tel que
$$div_S(f)=\sum\limits_{i=1}^n s_iC_i+\sum\limits_{s=1}^r m_sF_s+\sum\limits_{j=1}^t n_jD_j+lE'$$
où $D_1,\ldots D_t$ sont des courbes irréductibles, distinctes de $C_1,\ldots C_n, F_1,\ldots F_r$, $(n_j,l)=1$ et la spécialisation de $\alpha$ en $D_j$ est dans $H^2_{nr}(\kappa(D_j)/D_j, \Zl)$.\\}}
\rem{La spécialisation $\bar{\alpha}$ de $\alpha$ en $D_j$ est définie par $\bar{\alpha}=\partial_{\mathcal {O}_{S,D_j}}(\alpha\cup \pi)$ où $\pi$ est une uniformisante de $\mathcal{O}_{S,D_j}$, $\kappa(D_j)$ est son corps résiduel. Le groupe $H^2_{nr}(\kappa(D_j)/D_j, \Zl)$ désigne ici le sous-groupe de $H^2(\kappa(D_j), \Zl)$ formé des éléments qui sont non ramifiés pour toute valuation discrète de $\kappa(D_j)$ au-dessus d'un point fermé de $D_j$. Cette définition ne nécessite pas d'hypothèse de régularité de $D_j$. \\}

On passe maintenant à la preuve du principe local-global de \cite{PS}.\\

\textit{Démonstration du théorème \ref{lg}.}

L'implication $(i)\Rightarrow (ii)$ est immédiate. Montrons  que $(ii)\Rightarrow (i)$. On dispose alors des éléments $f_x$ dans le complété $K_x$ de $K$ pour tout  point $x\in S$ de codimension $1$. Notons $\kappa(x)$ le corps résiduel de $K_x$.

Soit $\mathcal C$ un ensemble fini de points de codimension $1$ de $S$, qui consiste exactement des points de $ram_S(\alpha)\cup ram_S(\xi)$. Par l'approximation faible, il existe un élément $f\in K^*$ tel que, pour tout $x\in \mathcal C$, $f/f_x$ est une puissance $l$-ième dans $K_x$.
On écrit :
$$div_S(f)=C-\sum\limits_{s=1}^r m_sF_s+lE$$
où $C$ est à support dans $\mathcal C$ et le support de $E$ est disjoint de $F_s$, $s=1,\ldots r$.

Ensuite, par l'approximation faible, on choisit $u\in K^*$ tel que
\begin{itemize}
\item[(i)] la valuation de $u$ en $F_i$ est $m_i$;
\item[(ii)] $u$ est une unité en chaque point $x\in \mathcal C$ et l'image $\bar u_x$ de $u$ dans $\kappa(x)^*/(\kappa(x)^*)^l$ est
$\bar u_x=\left\{
            \begin{array}{ll}
              \partial_x(\alpha), & x\in ram_S(\alpha); \\
              \hbox{non nul}, & \hbox{sinon.}
            \end{array}
          \right.
$
\end{itemize}

On pose $L=K(\sqrt[l]{u})$. On note $Y\stackrel{\pi}{\to} S$ la normalisation de $S_L$.  Notons que la condition $(m_j,l)=1$ assure que $Y\stackrel{\pi}{\to} S$ est ramifié en $F_i$. Ainsi on trouve une seule courbe $\tilde F_i$ dans la préimage de $F_i$ et, de plus, $\kappa(\tilde F_i)=\kappa(F_i)$. Soit $\tilde Y\stackrel{\eta}{\to} Y$ un modèle régulier de $Y$, tel que $ram_{\tilde Y}\alpha_L$ et les transformés stricts des  $\tilde F_i$, que l'on note aussi $\tilde F_i$,  soient à croisements normaux.

On applique ensuite la proposition \ref{div} à $\tilde Y$, $\alpha_L$ et $\tilde F_i$. On trouve $g\in L^*$ tel que
$$div_{\tilde Y}(g)=C'+\sum\limits_{s=1}^r m_s\tilde F_s+\sum\limits_{j=1}^t n_jD_j+lE'$$
où $C' $ est à support dans $ram_{\tilde Y}\alpha_L$. De plus, la spécialisation de $\alpha_L$ en $D_j$ est dans $H^2_{nr}(\kappa(D_j)/D_j, \Zl)$.

Montrons que $\xi'\eqdef\xi-\alpha\cup fN_{L/K}(g)$ convient : $\xi'\in\Hn(K/S,\Zl)$.
Soit $x\in S^{(1)}$. On a trois cas à considérer :
\begin{enumerate}
\item si $x\notin  ram_S(\alpha)\cup ram_S(\xi)\cup Supp (fN_{L/K}(g))$, $\xi'$ est non ramifiée en $x$.
\item Supposons que $x\in  ram_S(\alpha)\cup ram_S(\xi)$. D'après le choix de $f$, on a  $\partial_x(\xi-\alpha\cup f)=\partial_x(\xi-\alpha\cup f_x)=0$. Montrons que $\partial_x(\alpha\cup N_{L/K}(g))=0$. On étend la valuation sur $K$ donnée par $x$ en une valuation $v$ sur $L$. D'après le choix de $u$, $\alpha_L$ est triviale sur $L_v$. Ainsi  $\partial_x(\alpha\cup N_{L/K}(g))=\partial_x(Cores_{L_v/K_x}(\alpha_L\cup g))=0$.
\item Supposons que $x\in Supp (fN_{L/K}(g))\setminus(ram_S(\alpha)\cup ram_S(\xi))$. On a alors $\partial_x(\xi')=\partial_x(\alpha\cup fN_{L/K}(g))$. On a \\ $div_S(fN_{L/K}(g))=div_S(f)-\pi_*\eta_*(div_{\tilde Y}(g))$.
Puisque $\kappa(\tilde F_i)=\kappa(F_i)$, on en déduit $$div_S(fN_{L/K}(g))=C''+\sum\limits_{j=1}^t n_j\pi_*\eta_*(D_j)+lE'',$$ où $C''$ est supporté sur $\mathcal C$.

Soit $D_j'=\pi(\eta(D_j))$. Si $\pi_*\eta_*(D_j)=cD_j'$ est non nul et si $l\nmid c$, on a nécessairement que $\kappa(D_j')=\kappa(D_j)$. Ainsi soit $\partial_x(\alpha\cup fN_{L/K}(g))$ est nul, soit c'est un multiple de la spécialisation de $\alpha$ en $D_j'$, qui est dans $H^2_{nr}(\kappa(D_j)/D_j, \Zl)$ par le choix de $g$. Puisque $D_j$ est une courbe sur un corps fini,    $H^2_{nr}(\kappa(D_j)/D_j, \Zl)=0$. On obtient donc $\partial_x(\xi')=0$, ce qui termine la preuve. \qed
\end{enumerate}

\small{

}

\begin{thebibliography}{25}
\bibitem[A]{A} M. Artin, \emph{Left ideals in maximal orders},  Brauer groups in ring theory and algebraic geometry (Wilrijk, 1981),  pp. 182--193, Lecture Notes in Math., \textbf{917}, Springer, Berlin-New York, 1982.
\bibitem[CT]{CT} J.-L. Colliot-Th\'el\`ene, \emph{Birational invariants, purity and the Gersten conjecture},  $K$-theory and algebraic geometry: connections with quadratic forms and division algebras (Santa Barbara, CA, 1992),  1--64,Proc. Sympos. Pure Math., \textbf{58}, Part 1, Amer. Math. Soc., Providence, RI, 1995.
\bibitem[CTOj]{CTOj} J.-L. Colliot-Th\'el\`ene et M. Ojanguren, \emph{Vari\'et\'es unirationnelles non rationnelles: au-del\`a de l'exemple d'Artin et Mumford},
Invent. Math. \textbf{97} (1989), no. 1, 141--158.
\bibitem[CTSS]{CTSS} J.-L. Colliot-Th\'el\`ene, J.-J. Sansuc et C. Soulé, \emph{Torsion dans le groupe de Chow de codimension deux},  Duke Math. J.  \textbf{50}  (1983),  no. 3, 763--801.
Doc. Math. \textbf{1} (1996), No. 17, 395--416.
\bibitem[Fr]{Fr} E. Frossard, \emph{Fibres dégénérées des schémas de Severi-Brauer d'ordres},   J. Algebra  \textbf{198}  (1997),  no. 2, 362--387.
\bibitem[K1]{K1} B. Kahn, \emph{Motivic cohomology of smooth geometrically cellular varieties},  Algebraic $K$-theory (Seattle, WA, 1997),  149--174, Proc. Sympos. Pure Math., \textbf{67}, Amer. Math. Soc., Providence, RI, 1999.
\bibitem[K2]{K2} B. Kahn, \emph{Cohomological approaches to $SK_1$ and $SK_2$ of central simple algebras},
Documenta Mathematica,  Extra Volume: Andrei A. Suslin's Sixtieth Birthday (2010), 317--369.
\bibitem[MS]{MS} A. S. Merkur'ev and A. A. Suslin, \emph{$K$-cohomology of Severi-Brauer varieties and the norm residue homomorphism},
Izv. Akad. Nauk SSSR Ser. Mat. \textbf{46} (1982), no. 5, 1011--1046, 1135--1136.
\bibitem[PS]{PS} R. Parimala and V. Suresh, \emph{Degree three cohomology of function fields of surfaces}, 	arXiv:1012.5367v1 (2010).
\bibitem[Pe]{Pe} E. Peyre, \emph{Products of Severi-Brauer varieties and Galois cohomology},   $K$-theory and algebraic geometry: connections with quadratic forms and division algebras (Santa Barbara, CA, 1992),  369--401,
Proc. Sympos. Pure Math., \textbf{58}, Part 2, Amer. Math. Soc., Providence, RI, 1995.
\bibitem[S1]{S1}  D. J. Saltman, \emph{Cyclic algebras over $p$-adic curves},  J. Algebra  \textbf{314}  (2007),  no. 2, 817--843.
\bibitem[S2]{S2} D. J. Saltman,  \emph{Division algebras over surfaces},  J. Algebra  \textbf{320}  (2008),  no. 4, 1543--1585.
\bibitem[Se]{Se} J-P. Serre, \emph{Corps locaux},  Publications de l'Universit\'e de Nancago, No. VIII. Hermann, Paris, 1968.
\bibitem[Su]{Su} A. A. Suslin, 
\emph{Quaternion homomorphism for the field of functions on a conic}, 
Dokl. Akad. Nauk SSSR \textbf{265} (1982), no. 2, 292--296. 
\end{thebibliography}
\end{document}